\documentclass[12pt]{article}
\usepackage{amsmath}
\usepackage{amssymb}
\usepackage{latexsym}
\usepackage{amsthm}
\usepackage{tabularx}
\usepackage{booktabs}
\usepackage{mathrsfs}
\usepackage{enumitem}
\usepackage{ytableau}
\usepackage{graphicx}
\usepackage{algorithm,algorithmic}
\usepackage[colorlinks,
            linkcolor=red,
            anchorcolor=blue,
            citecolor=green]{hyperref}

\usepackage{epic}

\newtheorem{thm}{Theorem}[section]
\newtheorem{lem}[thm]{Lemma}

\newtheorem{cor}[thm]{Corollary}

\allowdisplaybreaks

\begin{document}
\begin{center}
{\large \bf  The equivariant inverse Kazhdan-Lusztig polynomials of uniform matroids}
\end{center}

\begin{center}
Alice L.L. Gao$^1$, Matthew H.Y. Xie$^2$ and Arthur L.B. Yang$^{3}$\\[6pt]

School of Mathematics and Statistics,\\
Northwestern Polytechnical University, Xi'an, Shaanxi 710072, P.R. China

College of Science, \\
Tianjin University of Technology, Tianjin 300384, P. R. China

Center for Combinatorics, LPMC\\
Nankai University, Tianjin 300071, P. R. China\\[6pt]

Email: $^{1}${\tt llgao@nwpu.edu.cn},
      $^{2}${\tt xie@email.tjut.edu.cn},
     $^{3}${\tt yang@nankai.edu.cn}
\end{center}
\noindent\textbf{Abstract.}
Motivated by the concepts of the inverse Kazhdan-Lusztig polynomial and the equivariant Kazhdan-Lusztig polynomial,
Proudfoot defined the equivariant inverse Kazhdan-Lusztig polynomial for a matroid. In this paper, we show that the
equivariant inverse Kazhdan-Lusztig polynomial of a matroid is very useful for determining its equivariant Kazhdan-Lusztig polynomials, and we determine the equivariant inverse Kazhdan-Lusztig polynomials for Boolean matroids and  uniform matroids.
As an application, we give a new proof of Gedeon, Proudfoot and Young's formula for the equivariant  Kazhdan-Lusztig polynomials of uniform matroids. Inspired by Lee, Nasr and Radcliffe's combinatorial interpretation for the ordinary
Kazhdan-Lusztig polynomials of uniform matroids, we further present a new formula for the corresponding equivariant Kazhdan-Lusztig polynomials.

\noindent \emph{AMS Classification 2020:} {05B35, 05E05, 20C30} 




\noindent \emph{Corresponding Author:} Matthew H.Y. Xie, xie@email.tjut.edu.cn

\section{Introduction}

Since the introduction of the matroid Kazhdan-Lusztig polynomials, due to Elias, Proudfoot and Wakefield \cite{ElisKL2016adv},
these polynomials have attracted much attention, for instance see \cite{Gedeon2017ejc, Gedeon_proud_2017sem, lxy_2018arxiv, glxyz_2018arxiv, BV_2020EJC, W_2018EJC} and references therein.
As noted by Proudfoot \cite{ proudfoot2017arxiv_alge-geo}, the Kazhdan-Lusztig polynomials of matroids can also be considered  as a special case of the Kazhdan-Lusztig-Stanley polynomials, which were first introduced by Stanley \cite{stanley1992kls} and further studied by Brenti \cite{brenti1999kls1, brenti1999kls2}. Within the framework of the Kazhdan-Lusztig-Stanley theory, Gao and Xie \cite{GaoXie_2020_inv} defined the inverse Kazhdan-Lusztig polynomial of a matroid.
To study the properties of the matroid Kazhdan-Lusztig polynomials, Gedeon, Proudfoot, and Young \cite{Geden2017jcta} introduced the concept of the equivariant Kazhdan-Lusztig polynomials of matroids. In \cite{P_2020} Proudfoot further developed the equivariant Kazhdan-Lusztig-Stanley theory, noted that the equivariant Kazhdan-Lusztig polynomial of a matroid can be realized as an equivariant Kazhdan-Lusztig-Stanley polynomial, and further defined the equivariant inverse Kazhdan-Lusztig polynomials of a matroid.
Gao and Xie \cite{GaoXie_2020_inv} showed that the inverse Kazhdan-Lusztig polynomials of uniform matroids are not only easy to compute, but also very helpful for determining their Kazhdan-Lusztig polynomials. The main objective of this paper is to develop their idea for the equivariant case of uniform matroids. To this end, we compute the equivariant inverse Kazhdan-Lusztig polynomials for Boolean matroids and uniform matroids, and further use them to give an alternative proof of Gedeon, Proudfoot and Young's formula  \cite{Geden2017jcta} for the equivariant Kazhdan-Lusztig polynomials of uniform matroids.

Let us first follow Gedeon, Proudfoot, and Young \cite{Geden2017jcta} to recall some related definitions and notations on equivariant matroid Kazhdan-Lusztig polynomials.
Let $M$ be a matroid  on the ground set $\mathcal{I}$,  and let $W$ be a finite group acting on $\mathcal{I}$ and preserving $M$. Gedeon, Proudfoot, and Young referred to this collection of data as an equivariant matroid $W \curvearrowright M$.
Let $\mathrm{VRep}(W)$ denote the ring of virtual representations of $W$ over $\mathbb{Z}$ with coefficients in $\mathbb{Z}$, and let $\mathrm{VRep}(W)[t]$ denote the polynomial ring in $t$ over $\mathrm{VRep}(W)$.
 Gedeon, Proudfoot and Young defined the equivariant characteristic polynomial of  $W \curvearrowright M$ as
\begin{align}\label{equi-charac-defi}
H_M^W(t):=\sum_{i=0}^{\mathrm{rk} \, M}(-1)^it^{\mathrm{rk} \, M-i} OS^W_{M,i} \in \mathrm{VRep}(W)[t],
\end{align}
where  $OS^W_{M,i}$ is the degree $i$ part of the Orlik-Solomon algebra of $M$, a natural representation of $W$ induced by its action on $M$. Note that the graded dimension of
$H_M^W(t)$  is just the usual characteristic polynomial $\chi_M(t)\in \mathbb{Z}[t]$.

Let $L(M)$ denote the lattice of flats of $M$. For any $F\in L(M)$, denote the localization of $M$ at $F$ by
$M_F$, the matroid on the ground set $F$ whose lattice of flats is isomorphic to $L_F:=\{G \in L ~|~ G \leq F\}.$
Dually, denote the contraction of $M$ at $F$ by $M^F$, the matroid on the ground set $\mathcal{I} \backslash F$ with its lattice of flats isomorphic to $L^F:=\{G \in L ~|~ G \geq F\}.$
Moreover, let $W_F \subset W$ be the stablizer of $F$ and let $\mathrm{rk}~M$ denote the rank of $M$.
Gedeon, Proudfoot, and Young \cite[Theorem 2.8]{Geden2017jcta} showed that
there is a unique way to assign to each equivariant matroid $W\curvearrowright M$ an element $P_M^W(t) \in \mathrm{VRep(W)}[t]$, called the equivariant Kazhdan-Lusztig polynomial, such that the following  conditions are satisfied:
\begin{itemize}
\item[(1).] If $\mathrm{rk}~M=0$, then $\mathrm{deg}~P_M^W(t)=0$, and $P_M^W(t)$ is the trivial representation of $W$.
\item[(2).] If $\mathrm{rk}~M>0$, then $\mathrm{deg}~P_M^W(t)<\frac{1}{2}~ \mathrm{rk}~M$.
\item[(3).] For every $M$,
\begin{align}\label{defi-KL_poly-matroid_Proud}
t^{\mathrm{rk}~M}    P_M^W(t^{-1})=\sum_{[F]\in L(M)/W}
\mathrm{Ind}_{W_F}^W \Big( H_{M_F}^{W_F}(t) \otimes P_{M^F}^{W_F}(t) \Big),
\end{align}	
where $L(M)/W$ denotes the set of orbits of the natural action of $W$ on $L(M)$.
\end{itemize}

The equivariant Kazhdan-Lusztig polynomials have been computed for uniform matroids \cite{Geden2017jcta},  q-niform matroids \cite{Proud2017arxiv-eq}, and  thagomizer matroids \cite{Gedeon2017ejc, xiezang2019pams}.
Gedeon, Proudfoot, and Young  \cite[Conjecture 2.13]{Geden2017jcta}  conjectured that for any equivariant matroid $W\curvearrowright M$ the coefficients of $P_M^W(t)$ are isomorphism classes of honest representations of $W$.
Recently, Braden, Huh,  Matherne,  Proudfoot and Wang \cite{BHMPW_2020} proved this conjecture by using the singular hodge theory.

As remarked by Gedeon, Proudfoot, and Young  \cite{Geden2017jcta}, their computation of the equivariant Kazhdan-Lusztig polynomials for uniform matroids relies on a guess on the generating function of these polynomials. Our proof given here is more direct, and uses the concept of the equivariant inverse Kazhdan-Lusztig polynomials of matroids. As noted by Proudfoot \cite[Section 4]{P_2020},
there is a unique way to assign to each equivariant matroid $W\curvearrowright M$ an element $Q_M^W(t) \in \mathrm{VRep}(W)[t]$ such that the following conditions are satisfied:
	\begin{itemize}
	\item[(a).] If $\mathrm{rk}\, M = 0$, then  $Q_M^W(t)$ is the trivial representation in degree $0$.
	\item[(b).] If $\mathrm{rk} \, M> 0$, then $\mathrm{deg}\,Q_M^W(t) < \frac{1}{2} \, \mathrm{rk}\, M  $.
	\item[(c).] For every $W\curvearrowright M$,
	\begin{align}\label{matroid-inv-Q-defi}
		t^{\mathrm{rk} \, M} \cdot  (-1)^{\mathrm{rk} \, M} Q^W_M(t^{-1})
		= \sum_{[F]\in L(M)/W}\mathrm{Ind}_{W_F}^{W}
 \Big(  (-1)^{\mathrm{rk} \, M_F}  Q_{M_F}^{W_F}(t)    \otimes   t ^{\mathrm{rk} \, M^F}H_{M^F}^{W_F}(t^{-1})\Big).
		\end{align}
	\end{itemize}
We call $Q_M^W(t)$ the equivariant inverse Kazhdan-Lusztig polynomial of  $W\curvearrowright M$.
Note that \eqref{matroid-inv-Q-defi} is equivalent to
		\begin{align}\label{inv_KL_equa_2_2}
(-1)^{\mathrm{rk} \, M}Q^W_M(t)
=\sum_{[F]\in L(M)/W}\mathrm{Ind}_{W_F }^{W}
 \Big((- t) ^{\mathrm{rk} \, M_F}Q_{M_F}^{W_F}(t^{-1})    \otimes  H_{M^F}^{W^F}(t) \Big);
\end{align}
By using the equivariant  Kazhdan-Lusztig-Stanley theory, it is easy to deduce that
 \begin{align}\label{KL-inverse-KL-or}
 \sum_{[F]\in L(M)/W} (-1)^{\mathrm{rk} \, M_F}
 \mathrm{Ind}_{W_F}^{W}
	\Big(  P_{M_F}^{W_F}(t)    \otimes  Q_{M^F}^{W_F}(t) \Big)
=0.
\end{align}
We would like to point out that \eqref{KL-inverse-KL-or} was used by Braden, Huh,  Matherne,  Proudfoot and Wang to define $Q_M^W(t)$.
Moreover, they showed that the coefficients of $Q_M^W(t)$ are isomorphism classes of honest representations of $W$. For notational convenience, let
$$\Hat{Q}_{M}^{W}(t)=(-1)^{\mathrm{rk} \, M}Q_{M}^{W}(t).$$
 Then the relations \eqref{inv_KL_equa_2_2} and \eqref{KL-inverse-KL-or} can be written as
\begin{align}\label{matroid-inv-Q-defi-hat}
\hat{Q}_M^W(t)
		=\sum_{[F]\in L(M)/W}\mathrm{Ind}_{W_F}^{W}  \Big( t ^{\mathrm{rk} \, F}\hat{Q}_{M_F}^{W_F}(t^{-1})    \otimes  H_{M^F}^{W_F}(t)\Big) ,
\end{align}
and
 \begin{align}\label{KL-inverse-KL}
 \sum_{[F]\in L(M)/W}\mathrm{Ind}_{W_F}^{W}
	\Big(  P_{M_F}^{W_F}(t)    \otimes   \hat{Q}_{M^F}^{W_F}(t) \Big)
=0.
\end{align}

This paper is organized as follows. In  Section \ref{sec-4-sym-func},  we will review some basic definitions and notations of the symmetric functions and give some results which will be used later. In Section  \ref{sec-5-uniform}, we first compute the equivariant inverse Kazhdan-Lusztig polynomials for Boolean matroids, and then compute these polynomials for uniform matroids, based on which we
present a new proof for Gedeon, Proudfoot, and Young's formula for the equivariant Kazhdan-Lusztig polynomial for uniform matroids, as well as a new formula of these equivariant Kazhdan-Lusztig polynomials.

\section{Symmetric functions}\label{sec-4-sym-func}
The aim of this section is to  review some basic definitions and results on symmetric functions. We refer the reader to Stanley \cite{stanley-enumer_com_book_2_1999}, Macdonald \cite{Macdonald_book_symme_2008} and Haglund \cite{Haglund_book_2008}  for undefined terminology from the theory of symmetric functions. We also prove some symmetric function identities which will be used in the evaluation of the (inverse) equivariant Kazhdan-Lusztig polynomials for uniform matroids.

\subsection{The Frobenius characteristic map}

Let $\Lambda_n$ denote the $\mathbb{Z}$-module of symmetric functions of degree $n$ in the variables
$\mathrm{\bold x}=(x_1,x_2,\ldots)$. From the theory of symmetric functions,
the Frobenius characteristic map is an isomorphism
$$\mathrm{ch}: \mathrm{VRep}(S_n)\longrightarrow \Lambda_n,$$
which maps the irreducible representation $V_{\lambda}$ of $S_n$ to the Schur function $s_{\lambda}(\mathrm{\bold x})$ for each partition $\lambda$ of $n$. In particular, the image of the trivial representation $V_{(n)}$ is $s_{(n)}(\mathrm{\bold x})=h_n(\mathrm{\bold x})$, which is called the complete symmetric function, and the image of the representation $V_{(1^n)}$ is $s_{(1^n)}(\mathrm{\bold x})=e_n(\mathrm{\bold x})$, which is called the elementary symmetric function.
Moreover, the image of the skew  Specht module $V_{\lambda/\mu}$ under $\mathrm{ch}$ is the skew Schur function $s_{\lambda/\mu}(\mathrm{\bold x}).$
The definition of $\mathrm{ch}$ carries over directly from $\mathrm{VRep}(S_n)$ to $\mathrm{VRep}(S_n)[t]$.
It has the property that, given two graded virtual representations $V_1 \in \mathrm{VRep}(S_{n_1})[t]$ and $V_2 \in \mathrm{VRep}(S_{n_2})[t]$, we have
\begin{align}\label{eq-preserving-ch}
\mathrm{ch}\, \mathrm{Ind}_{S_{n_1}\times S_{n_2}}^{S_{n_1+n_2}}(V_1\otimes V_2)=\mathrm{ch}(V_1)\mathrm{ch}(V_2).
\end{align}

\subsection{Plethystic}

Here we adopt the notation of Haglund \cite{Haglund_book_2008}. Let $E=(t_1,t_2,t_3,\ldots)$ be a formal series of rational functions in the parameters. Let $p_k[E]$ denote the plethystic substitution of $E$ into the $k$-th power sum $p_k$, i.e.,
$$p_k[E]=E(t_1^k,t_2^k,\ldots).$$
For any symmetric function $f$, suppose that $f=\sum_{\lambda}c_{\lambda}p_{\lambda}=\sum_{\lambda}c_{\lambda}\prod_{i}p_{\lambda_i}$, and then define
$$f[E]=\sum_{\lambda}c_{\lambda}\prod_{i}p_{\lambda_i}[E].$$
Note that if $X=\sum_ix_i$, then $p_k[X]=p_k(\mathrm{\bold x})$, from which we get
 $f[X]=f(\mathrm{\bold x})$ for any symmetric function $f$.
It is also easy to see that  $p_k[-X]=-p_k(\mathrm{\bold x})=-p_k[X]$ and $e_m[tX]=t^{m}e_{m}(\mathrm{\bold x})$.


\begin{lem}[\cite{Loehr_plethystic},Section 3.3]
\label{A-B}
Let $E=E(t_1,t_2,t_3,\ldots)$ and $F=F(w_1,w_2,w_3,\ldots)$ be two formal series  of rational functions in the parameters $t_1,t_2,\ldots$ and $w_1,w_2,\ldots$. Then
$$e_m[E-F]=\sum_{j=0}^m(-1)^{m-j}e_j[E]h_{m-j}[F].$$
\end{lem}

\subsection{The Pieri rule on Schur functions}

In this subsection we present three symmetric function identities, which might be known.
To be self-contained here, we will use the Pieri rule to give proofs of these identities.

Let us first recall the well known Pieri rule on Schur functions; see Stanley \cite{stanley-enumer_com_book_2_1999} for details.
For any $i\geq 1$, we have
\begin{align*}
s_{(i)}(\mathrm{\bold x})s_{\lambda}(\mathrm{\bold x})=\sum_{\mu}s_{\mu}(\mathrm{\bold x})
\end{align*}
 summing over all partitions $\mu$ such that $\mu/\lambda$ is a horizontal strip of size $i$.  Meanwhile, we have
 \begin{align*}
s_{(1^i)}(\mathrm{\bold x})s_{\lambda}(\mathrm{\bold x})=\sum_{\mu}s_{\mu}(\mathrm{\bold x})
\end{align*}
 summing over all partitions $\mu$ such that $\mu/\lambda$ is a vertical strip of size $i$.

We proceed to prove the first symmetric function identity, which is stated as follows.
\begin{lem} \label{schur-piere-sub-2}
For $m\geq2$, $i \geq 0$   and $j -i\geq 2$, we have
\begin{equation}\label{schur-piere-sub-2-eq}
s_{(1^{i+1})}(\mathrm{\bold x})s_{(m,1^{j-1})}(\mathrm{\bold x})-s_{(1^i)}(\mathrm{\bold x})s_{(m,1^j)}(\mathrm{\bold x})
=s_{(m+1,2^i,1^{j-i-1})}(\mathrm{\bold x})+s_{(m,2^{i+1},1^{j-i-2})}(\mathrm{\bold x}).
\end{equation}
\end{lem}

 \proof
For  $m\geq 2$ and $i,j \geq 0$,  by Pieri rule, we have
\begin{align*}
	s_{(1^i)}(\mathrm{\bold x})s_{(m,1^j)}(\mathrm{\bold x})=\sum_{x=0}^{\mathrm{min}\, \{i-1,j\}}s_{(m+1,2^x,1^{i+j-1-2x})}(\mathrm{\bold x})+\sum_{x=0}^{\mathrm{min}\, \{i,j\}}s_{(m,2^x,1^{i+j-2x})}(\mathrm{\bold x}).
\end{align*}
It follows that, for $m\geq 2$, $i\geq 0$ and $ j \geq 1$,
\begin{align*}
&s_{(1^{i+1})}(\mathrm{\bold x})s_{(m,1^{j-1})}(\mathrm{\bold x})-
s_{(1^i)}(\mathrm{\bold x})s_{(m,1^j)}(\mathrm{\bold x})\\
&~~~~~=\sum_{x=0}^{\min\{i,j-1\}}s_{(m+1,2^x,1^{i+j-1-2x})}(\mathrm{\bold x})+\sum_{x=0}^{\min\{i+1,j-1\}}s_{(m,2^x,1^{i+j-2x})}(\mathrm{\bold x})\\
&~~~~~~~~~~~~~~-\sum_{x=0}^{\min\{i-1,j\}}s_{(m+1,2^x,1^{i+j-1-2x})}(\mathrm{\bold x})-\sum_{x=0}^{\min\{i,j\}}s_{(m,2^x,1^{i+j-2x})}(\mathrm{\bold x}).
\end{align*}
In view of that $j-i\geq 2$, and hence $i\leq j-2$, the above four summations reduce to the right hand side of \eqref{schur-piere-sub-2-eq}. This completes the proof.
\qed

The second symmetric function identity we are to prove is as follows.

\begin{lem}\label{schur-pieri-3-lem}
For $n\geq 0 $, $m\geq 2$ and $j\geq 1$, we have
$$\sum_{i=0}^n(-1)^is_{(i)}(\mathrm{\bold x})s_{(m,2^j,1^{n-i})}(\mathrm{\bold x})=(-1)^n\sum_{x=0}^{\mathrm{min}\, \{m-2,n \}}s_{(m+n-x,2+x,2^{j-1})}(\mathrm{\bold x}).$$
\end{lem}

\proof
We may assume that $n\geq 1$ since the case of $n=0$ is obvious. Given $m\geq 2$, $j\geq 1$ and $n\geq 1 $, let $\lambda^i=(m,2^j,1^{n-i})$ for $0\leq i \leq n$.
According to Pieri rule, we have
\begin{align}\label{Schur-piere-h_i}
s_{(i)}(\mathrm{\bold x})s_{(m,2^j,1^{n-i})}(\mathrm{\bold x})=\sum_{\mu}s_{\mu}(\mathrm{\bold x}),
\end{align}
where the summation ranges over all partitions $\mu \vdash m+2j+n$ such that $\mu / \lambda^i$ is a horizontal strip of size $i$.
Considering the shape of $\lambda^i$, such a partition $\mu$ must be of the form $(A,B,2^u,1^{v})$ for some
$A\geq m,\, B \geq 2$, $u=j-1$ or $j$, and nonnegative $v \geq n-i-1$. Assume that
$$\sum_{i=0}^n(-1)^is_{(i)}(\mathrm{\bold x})s_{(m,2^j,1^{n-i})}(\mathrm{\bold x})=\sum_{\mu}c_{\mu}s_{\mu}(\mathrm{\bold x}),$$
where $\mu=(A,B,2^u,1^{v})\vdash m+2j+n$ for $A\geq m,\, B \geq 2$, $u=j-1$ or $j$, and $v \geq 0$.
It remains to determine the coefficient $c_{\mu}$.

Now fix a partition $\mu=(A,B,2^u,1^{v})$. It is easy to see that
$$c_{\mu}=\sum_{i}(-1)^i,$$
where the sum is over all $i$ such that $s_{(i)}(\mathrm{\bold x})s_{(m,2^j,1^{n-i})}(\mathrm{\bold x})$ contains $s_{\mu}(\mathrm{\bold x})$ in \eqref{Schur-piere-h_i}. We claim that $c_{\mu}=0$ for $v\geq 1$. There are two cases to consider.

\begin{itemize}
\item[(i)]
Suppose that $v\geq 1$ and $u=j-1$. To guarantee that $s_{\mu}(\mathrm{\bold x})$ appears in the Schur expansion of $s_{(i)}(\mathrm{\bold x})s_{\lambda^i}(\mathrm{\bold x})$, there are only two choices for $i$: $i=(A-m)+(B-2)$ or $i=(A-m)+(B-2)+1$.
Precisely, we have $\lambda^{i}=(m,2,2^{j-1}, 1^v)$ or $\lambda^{i}=(m,2,2^{j-1}, 1^{v-1})$. Thus, $c_{\mu}=0.$

\item[(ii)] Suppose that $v\geq 1$ and $u=j$. To guarantee that $s_{\mu}(\mathrm{\bold x})$ appears in the Schur expansion of $s_{(i)}(\mathrm{\bold x})s_{\lambda^i}(\mathrm{\bold x})$, there are only two choices for $i$: $i=(A-m)+(B-2)+1$ or $i=(A-m)+(B-2)+2$.
Precisely, we have $\lambda^{i}=(m,2,2^{j-1}, 1^{v+1})$ or $\lambda^{i}=(m,2,2^{j-1}, 1^{v})$. Thus, $c_{\mu}=0.$
\end{itemize}

We proceed to determine the coefficient $c_\mu$ for $\mu$ being of the form $(A,B,2^u)$. We would like to point out that for $0 \leq i \leq n-2$ the product $s_{(i)}(\mathrm{\bold x})s_{(m,2^j,1^{n-i})}(\mathrm{\bold x})$ does not contribute such $s_{\mu}(\mathrm{\bold x})$
due to the fact that the number of occurrences of $1$ in  partition $(m,2^j,1^{n-i})$ is $n-i \geq 2$, 
which implies that, for each $s_{\mu}(\mathrm{\bold x})$ in the Schur expansion of the product $s_{(i)}(\mathrm{\bold x})s_{(m,2^j,1^{n-i})}(\mathrm{\bold x})$, the partition $\mu$ must contain $1$ as a part by the Pieri rule. It remains to consider the contribution of $s_{(i)}(\mathrm{\bold x})s_{(m,2^j,1^{n-i})}(\mathrm{\bold x})$ when $i=n-1$ and $i=n$.
Note that, by the Pieri rule,
\begin{align*}
s_{(n-1)}(\mathrm{\bold x})s_{(m,2^j,1)}(\mathrm{\bold x})
=\sum_{x=0}^{\mathrm{min}\, \{m-2,n-2 \}}s_{(m+n-2-x,2+x,2^j)}(\mathrm{\bold x})+\sum_{\rho}s_{\rho}(\mathrm{\bold x})
\end{align*}
and
\begin{align*}
s_{(n)}(\mathrm{\bold x})s_{(m,2^j)}(\mathrm{\bold x})
=&\sum_{x=0}^{\mathrm{min}\, \{m-2,n \}}s_{(m+n-x,2+x,2^{j-1})}(\mathrm{\bold x})
+\sum_{x=0}^{\mathrm{min}\, \{m-2,n-2 \}}s_{(m+n-2-x,2+x,2^{j})}(\mathrm{\bold x})\\
&~~~~~~~~~+\sum_{\tau}s_{\tau}(\mathrm{\bold x}),
\end{align*}
where $\rho$ and $\tau$ are some partitions each of which contains $1$ as a part.
Hence, we have
\begin{align*}
\sum_{i=0}^n(-1)^i&s_{(i)}(\mathrm{\bold x})s_{(m,2^j,1^{n-i})}(\mathrm{\bold x})\\
=&(-1)^{n-1}\sum_{x=0}^{\mathrm{min}\, \{m-2,n-2 \}}s_{(m+n-2-x,2+x,2^j)}(\mathrm{\bold x})\\
&~~+(-1)^n \Big(\sum_{x=0}^{\mathrm{min}\, \{m-2,n \}}s_{(m+n-x,2+x,2^{j-1})}(\mathrm{\bold x})
+\sum_{x=0}^{\mathrm{min}\, \{m-2,n-2 \}}s_{(m+n-2-x,2+x,2^{j})}(\mathrm{\bold x}) \Big)\\
=&(-1)^n\sum_{x=0}^{\mathrm{min}\, \{m-2,n \}}s_{(m+n-x,2+x,2^{j-1})}(\mathrm{\bold x}).
\end{align*}
This completes the proof.
\qed

The third symmetric function identity we are to prove is as follows.
\begin{lem} \label{pieri-4-lem-equi-KL-uniform}
For $n \geq 0$ and $m\geq 1$, we have
$$\sum_{i=0}^{n}(-1)^is_{(i)}(\mathrm{\bold x}) s_{(m+1,1^{n-i})}(\mathrm{\bold x})=(-1)^ns_{(m+n+1)}(\mathrm{\bold x}).$$
\end{lem}
\proof
We may assume that $n\geq 1$ since the case of $n=0$ is obvious.
Now apply the Pieri rule to $s_{(i)} (\mathrm{\bold x})s_{(m+1,1^{n-i})}(\mathrm{\bold x})$.
For $i=n$, we have
 $$s_{(n)}(\mathrm{\bold x}) s_{(m+1)}(\mathrm{\bold x})=\sum_{x=0}^{\mathrm{min}\, \{m+1, n \}}s_{(m+1+n-x,x)}(\mathrm{\bold x}).$$
For $1 \leq i \leq n-1$, there holds 
\begin{align*}
s_{(i)} (\mathrm{\bold x})s_{(m+1,1^{n-i})}(\mathrm{\bold x})= \sum_{x=0}^{\mathrm{min}\,\{m, i-1\}}s_{(m+i-x, x+1, 1^{n-i})}(\mathrm{\bold x})+\sum_{x=0}^{\mathrm{min}\, \{m, i \}}s_{(m+1+i-x, x+1, 1^{n-i-1})}(\mathrm{\bold x}),
\end{align*}
and hence
\begin{align*}
\sum_{i=1}^{n-1}&(-1)^i s_{(i)} (\mathrm{\bold x})s_{(m+1,1^{n-i})}(\mathrm{\bold x})\\[5pt]
= &\sum_{i=1}^{n-1}(-1)^i\sum_{x=0}^{\mathrm{min}\,\{m, i-1\}}s_{(m+i-x, x+1, 1^{n-i})}(\mathrm{\bold x})+\sum_{i=1}^{n-1}(-1)^i \sum_{x=0}^{\mathrm{min}\, \{m, i \}}s_{(m+1+i-x, x+1, 1^{n-i-1})}(\mathrm{\bold x})\\[5pt]
= &-\sum_{i=0}^{n-2}(-1)^i\sum_{x=0}^{\mathrm{min}\,\{m, i\}}s_{(m+1+i-x, x+1, 1^{n-i-1})}(\mathrm{\bold x})+\sum_{i=1}^{n-1}(-1)^i \sum_{x=0}^{\mathrm{min}\, \{m, i \}}s_{(m+1+i-x, x+1, 1^{n-i-1})}(\mathrm{\bold x})\\[5pt]
=&-s_{(m+1, 1^{n})}(\mathrm{\bold x})+(-1)^{n-1}  \sum_{x=1}^{\mathrm{min}\,\{m+1, n \}}  s_{(m+n+1-x, x)}(\mathrm{\bold x}).
\end{align*}
Therefore, we obtain 
\begin{align*}
\sum_{i=0}^{n}&(-1)^is_{(i)}(\mathrm{\bold x}) s_{(m+1,1^{n-i})}(\mathrm{\bold x})\\[3pt]
=&(-1)^ns_{(n)}(\mathrm{\bold x}) s_{(m+1)}(\mathrm{\bold x})+s_{(m+1,1^{n})}(\mathrm{\bold x})+\sum_{i=1}^{n-1}(-1)^is_{(i)}(\mathrm{\bold x}) s_{(m+1,1^{n-i})}(\mathrm{\bold x})\\[5pt]
=&(-1)^n\sum_{x=0}^{\mathrm{min}\, \{m+1, n \}}s_{(m+1+n-x,x)}(\mathrm{\bold x})+(-1)^{n-1}  \sum_{x=1}^{\mathrm{min}\,\{m+1, n \}}  s_{(m+n+1-x, x)}(\mathrm{\bold x})\\[5pt]
=&(-1)^ns_{(m+n+1)}(\mathrm{\bold x}),
\end{align*}
as desired. This completes the proof. 
\qed

\subsection{The Littlewood-Richardson rule}

The Littlewood-Richardson rule gives a combinatorial interpretation of the Schur expansion of a skew Schur function.
There are many ways to state the Littlewood-Richardson rule; see Stanley \cite{stanley-enumer_com_book_2_1999} and references therein.
Here will use its Littlewood-Richardson tableaux version.

Recall that, given two partitions $\lambda$ and $\mu$ with $\mu \subseteq \lambda $ (i.e., $\mu_i \leq \lambda_i$ for all $i$), a semistandard Young tableau of shape $\lambda / \mu$ is defined to be an array $T=(T_{ij})$ of positive integers of shape $\lambda / \mu$ that is weakly increasing along every row and strictly increasing down every column. We say that $T$ has type $\alpha=(\alpha_1,\alpha_2,\ldots)$, if $T$ has $\alpha_i$ entries equal to $i$.
The reverse reading word of $T$ is the sequence of entries of $T$ obtained by concatenating the rows of $T$ from right to left, top to bottom. For example, the tableau
\begin{center}
\begin{ytableau}
  \none& \none&  \none&1& 1& 1& 2 \\
    \none&   \none& \none& 2& 2 \\
       \none&   \none& \none& 3& 3 \\
  1&2&2&4&4
\end{ytableau}
\end{center}
has the  reverse  reading word 2111~22~33~44221. 
We say that a word $a_1a_2\cdots a_n$ is a lattice permutation if in any initial factor $a_1a_2\cdots a_j$, the number of $i$'s is at least as great as the number of $i+1$'s (for all $i$). A Littlewood-Richardson tableau is a semistandard Young tableau $T$ such that its reverse reading word is a lattice permutation.

The well known Littlewood-Richardson rule is stated as follows.
\begin{thm}[{\cite[Section 7.10]{stanley-enumer_com_book_2_1999}}]\label{thm-lrrule}
If
$$s_{\lambda /\mu}(\mathrm{\bold x})=\sum_{\nu}c^{\lambda}_{\mu \nu}s_{\nu}(\mathrm{\bold x}),$$
then the Littlewood-Richardson coefficient $c^{\lambda}_{\mu \nu}$ is equal to the number of Littlewood-Richardson tableaux of shape $\lambda/ \mu$ and type $\nu$.
\end{thm}

\section{Uniform matroids}\label{sec-5-uniform}

Given a positive integer $d$ and a nonnegative integer $m$, let $U_{m,d}$ be the uniform matroid of rank $d$ on $m+d$ elements which admits an action of the symmetric group $S_{m+d}$. By using the generating functions and the Frobenius characteristic map, Gedeon,  Proudfoot, and Young \cite[Proposition 3.9]{ Geden2017jcta} obtained a formula for computing the equivariant Kazhdan-Lusztig polynomial for equivariant uniform matroids, which could be stated as follows.

\begin{thm}
[{\cite[Theorem 3.1]{Geden2017jcta}}]\label{main-eq-o}
	For any $m\geq 0$ and $d \geq 1$, we have
	\begin{align}\label{eq-ekl-uniform}
P_{U_{m,d}}^{S_{m+d}}(t)=V_{(m+d)}+\sum_{j=1}^{\lfloor (d-1)/2 \rfloor}t^j\sum_{x=1}^{\mathrm{min}\, \{m,d-2j \}}V_{(m+d-2j-x+1,x+1,2^{j-1})}.
\end{align}
\end{thm}

The main result of this section is as follows.

\begin{thm}\label{main-thm-inver-uniform}
For any equivariant uniform matroid $S_{m+d} \curvearrowright U_{m,d}$ with $m\geq 0$ and $d\geq 1$, we have
\begin{align}\label{formula-inverse-kl}
	Q_{U_{m,d}}^{S_{m+d}}(t)=\sum_{j=0}^{\lfloor (d-1)/2 \rfloor}V_{(m+1,2^{j},1^{d-2j-1})}t^j,
\end{align}
where $V_{(m+1,2^{j},1^{d-2j-1})}$ vanishes if $(m+1,2^{j},1^{d-2j-1})$ is not a valid partition.
\end{thm}

In subsection \ref{sec-4-Boolean} we first prove the above theorem for the case of $m=0$.
Then in subsection \ref{sec-5-uniform-inv} we prove Theorem \ref{main-thm-inver-uniform} for general $m$.
In subsection \ref{newproof} we use \eqref{formula-inverse-kl} to give a new proof of \eqref{eq-ekl-uniform}.
Finally in subsection \ref{newform} we present a new formula for $P_{U_{m,d}}^{S_{m+d}}(t)$, which refines
Lee, Nasr and Radcliffe's combinatorial interpretation for the ordinary Kazhdan-Lusztig polynomials of uniform matroids.

\subsection{Proof of Theorem \ref{main-thm-inver-uniform} for $m=0$}\label{sec-4-Boolean}

Given a positive integer $n$, let $B_n$ denote the Boolean matroid of rank $n$, which is equipped with a natural action of the symmetric group $S_{n}$. Gedeon,  Proudfoot, and Young \cite{Geden2017jcta} obtained a formula for computing the equivariant characteristic polynomials for
equivariant Boolean matroids, which could be written as follows in terms of the plethystic notation.

\begin{lem}[{\cite[Proposition 3.9]{Geden2017jcta}}]\label{prop-Boolean-Charac}
For any equivariant Boolean matroid $S_n\curvearrowright B_n$ with $n \geq 1$,
	we have
\begin{align}\label{eq-ecp-boolean}
\mathrm{ch} \, H_{B_n}^{S_n}(t)=h_n[(t-1)X].
\end{align}
\end{lem}

Note that the equivariant Boolean matroid
$S_n \curvearrowright B_n$ can be considered as the equivariant uniform matroid $S_n \curvearrowright U_{0,n}$.
Thus, the case $m=0$ of Theorem \ref{main-thm-inver-uniform} is equivalent to the following statement.

\begin{thm} \label{IKL_B_theo}
For any equivariant Boolean matroid $S_n\curvearrowright B_n$ with $n\geq 1$,
	we have $$Q_{B_n}^{S_n}(t)=V_{(1^n)}.$$
\end{thm}

\proof
It suffices to show that,  for $n\geq 1$,
$$\mathrm{ch} \, \Hat{Q}_{B_n}^{S_n}(t)=(-1)^ne_n(\mathrm{\bold x}).$$
Applying \eqref{matroid-inv-Q-defi-hat} to the equivariant Boolean matroid $S_n\curvearrowright B_n$, we obtain
\begin{align*}
\Hat{Q}_{B_n}^{S_n}(t)&=\sum_{[F]\in L(B_n)/S_n} \mathrm{Ind}_{(S_n)_F}^{S_n}
\Big( t^{\mathrm{rk} \, F}  \Hat{Q}_{(B_n)_F}^{(S_n)_F}(t^{-1})     \otimes  H_{(B_n)^F}^{(S_n)_F}(t)   \Big)\\
&=\sum_{i=0}^n\mathrm{Ind}_{S_i \times S_{n-i}}^{S_n} \Big( t^i \Hat{Q}_{B_i}^{S_i \times S_{n-i}}(t) \otimes H_{B_{n-i}}^{S_i \times S_{n-i}}(t)\Big),
\end{align*}
which is equivalent to
$$\Hat{Q}_{B_n}^{S_n}(t)=\sum_{i=0}^n\mathrm{Ind}_{S_i \times S_{n-i}}^{S_n} \Big( t^i \Hat{Q}_{B_i}^{S_i}(t) \otimes H_{B_{n-i}}^{S_{n-i}}(t)\Big).$$
From \eqref{eq-preserving-ch} and  \eqref{eq-ecp-boolean} it follows that
\begin{align*}
	\mathrm{ch} \, \Hat{Q}_{B_n}^{S_n}(t)=\sum_{i=0}^n t^i \mathrm{ch} \, \Hat{Q}_{B_i}^{S_i}(t^{-1}) \cdot \mathrm{ch} \, H_{B_{n-i}}^{S_{n-i}}(t)
	=\sum_{i=0}^n t^i \mathrm{ch} \, \Hat{Q}_{B_i}^{S_i}(t^{-1}) \cdot h_{n-i}[(t-1)X].
\end{align*}
Thus, we have
\begin{align}\label{Boolean-IKL-compute}
\mathrm{ch} \, \Hat{Q}_{B_n}^{S_n}(t)-t^n\mathrm{ch} \, \Hat{Q}_{B_n}^{S_n}(t^{-1})
=\sum_{i=0}^{n-1} t^i \mathrm{ch} \, \Hat{Q}_{B_i}^{S_i}(t^{-1}) \cdot h_{n-i}[(t-1)X].
\end{align}

Now we proceed to prove that $\mathrm{ch} \, \Hat{Q}_{B_n}^{S_n}(t)=(-1)^ne_n(\mathrm{\bold x})$ by induction on the value of $n$. For the base case, assume $n=1$.
Now \eqref{Boolean-IKL-compute} turns out to be
$$\mathrm{ch} \, \Hat{Q}_{B_1}^{S_1}(t)-t\mathrm{ch} \, \Hat{Q}_{B_1}^{S_1}(t^{-1})
=h_1[(t-1)X]=(t-1)h_1(\mathrm{\bold x}).$$
Since $\Hat{Q}_{B_1}^{S_1}(t)$ is of degree zero, we have
$$\mathrm{ch} \, \Hat{Q}_{B_1}^{S_1}(t)=-h_1(\mathrm{\bold x})=-e_1(\mathrm{\bold x}).$$
Assume the assertion for $n-1$, namely
\begin{align*}
\mathrm{ch} \,  \Hat{Q}_{B_i}^{S_i}(t)=(-1)^ie_i(\mathrm{\bold x})=(-1)^ie_i(X) \mbox{ for }0\leq i<n.
\end{align*}
From \eqref{Boolean-IKL-compute} we have
\begin{align*}
	 \mathrm{ch} \, \Hat{Q}_{B_n}^{S_n}(t)-t^n\mathrm{ch} \, \Hat{Q}_{B_n}^{S_n}(t^{-1})
	&=\sum_{i=0}^{n-1} t^i (-1)^i e_i(X) h_{n-i}[(t-1)X]\\
	&=\sum_{i=0}^{n} t^i (-1)^i e_i(X) h_{n-i}[(t-1)X]-t^n(-1)^ne_n[X].
	\end{align*}
Recall that   $t^i e_i[X]= e_i[tX]$ for $0 \leq i \leq n$, which tells us that
$$\mathrm{ch} \, \Hat{Q}_{B_n}^{S_n}(t)-t^n\mathrm{ch} \, \Hat{Q}_{B_n}^{S_n}(t^{-1})=(-1)^n\sum_{i=0}^{n}(-1)^{n-i} e_i[tX] h_{n-i}[(t-1)X]-t^n(-1)^ne_n[X].
$$
Taking $E=tX$ and $F=(t-1)X$ in Lemma \ref{A-B} leads to
\begin{align*}
	\mathrm{ch} \, \Hat{Q}_{B_n}^{S_n}(t)-t^n\mathrm{ch} \, \Hat{Q}_{B_n}^{S_n}(t^{-1})=(-1)^ne_n[X]-t^n(-1)^ne_n[X].
\end{align*}

In view of $ \mathrm{deg}\, \Hat{Q}_{B_n}^{S_n}(t)<\frac{n}{2}$, we find that
$$\mathrm{ch} \, \Hat{Q}_{B_n}^{S_n}(t)=(-1)^ne_n[X]=(-1)^ne_n(\mathrm{\bold x}),$$
as desired. This completes the proof.
\qed

\subsection{Proof of Theorem \ref{main-thm-inver-uniform} for general $m$}\label{sec-5-uniform-inv}

We proceed to prove Theorem \ref{main-thm-inver-uniform} for general $m$.

\noindent \textit{Proof of Theorem \ref{main-thm-inver-uniform}.}
It suffices to show for $0\leq j \leq \lfloor (d-1)/2 \rfloor$, the coefficient of $t^j$ in $ \mathrm{ch}~\Hat{Q}^{S_{m+d}}_{U_{m,d}}(t)$ is $$(-1)^ds_{(m+1,2^{j},1^{d-2j-1})}(\mathrm{\bold x}).$$
By applying \eqref{matroid-inv-Q-defi-hat} to the equivariant uniform matroid $S_{m+d}\curvearrowright U_{m,d}$, we obtain
\begin{align}\label{eq-qu}
\Hat{Q}_{U_{m,d}}^{S_{m+d}}(t)=\sum_{[F]\in L(U_{m,d})/S_{m+d}} \mathrm{Ind}_{(S_{m+d})_F}^{S_{m+d}}
\Big(  t^{\mathrm{rk}\, (U_{m,d})_F} \Hat{Q}_{(U_{m,d})_F}^{(S_{m+d})_F}(t^{-1})     \otimes  H_{(U_{m,d})^F}^{(S_{m+d})_F}(t)   \Big).
\end{align}
It is obvious that for any flat $F$ of $ L(U_{m,d})$ with  $\mathrm{rk} \, (U_{m,d})_F=i<d$ we have
\begin{align*}
(U_{m,d})_F \cong B_i,\qquad (U_{m,d})^F\cong U_{m,d-i},\qquad (S_{m+d})_F\cong S_i \times S_{m+d-i}.
\end{align*}
Note that there is only one flat of rank $d$ which is in fact
the ground set $\mathcal{I}$, for which there holds
\begin{align*}
(U_{m,d})_\mathcal{I}=U_{m,d},\qquad (U_{m,d})^\mathcal{I}\cong U_{0,0},\qquad (S_{m+d})_\mathcal{I}\cong S_{m+d},
\end{align*}
and hence $H_{(U_{m,d})^\mathcal{I}}^{(S_{m+d})_\mathcal{I}}(t)$ is just the trivial representation of $S_{m+d}$.
Thus, \eqref{eq-qu} can be rewritten as
\begin{align}\label{eq-qus}
\Hat{Q}_{U_{m,d}}^{S_{m+d}}(t)=\sum_{i=0}^{d-1} \mathrm{Ind}_{S_i \times S_{m+d-i}}^{S_{m+d}}
\Big( t^i \Hat{Q}_{B_i}^{S_i} (t^{-1}) \otimes H_{U_{m,d-i}} ^{S_{m+d-i}}(t) \Big)+ t^d \Hat{Q}_{U_{m,d}}^{S_{m+d}}(t^{-1}).
\end{align}
Applying the Frobenius characteristic map on both sides leads to
$$\mathrm{ch} \, \Hat{Q}^{S_{m+d}}_{U_{m,d}}(t)- t^d \mathrm{ch} \, \Hat{Q}^{S_{m+d}}_{U_{m,d}}(t^{-1})
=\sum_{i=0}^{d-1} t^i  \mathrm{ch} \, \Hat{Q}_{B_i}^{S_i}(t^{-1})\cdot \mathrm{ch} \, H_{U_{m,d-i}}^{S_{m+d-i}}(t).$$
By a result due to Gedeon,  Proudfoot, and Young \cite[Proposition 3.9 ]{ Geden2017jcta}, we know
\begin{align}
	\mathrm{ch} \, H_{U_{m,d-i}}^{S_{m+d-i}}(t)
	&=\sum_{j=0}^{d-i-1}(-1)^jt^{d-i-j} \Big(s_{(m+d-i-j,1^j)}(\mathrm{\bold x})+s_{(m+d-i-j+1,1^{j-1})}(\mathrm{\bold x})\Big)\nonumber\\
	&~~~~~~~~~~~~~+(-1)^{d-i}s_{(m+1,1^{d-i-1})}(\mathrm{\bold x}).\label{eq-hus}
	\end{align}
Meanwhile, by Theorem \ref{IKL_B_theo}, we have
\begin{align}\label{eq-qbs}
\mathrm{ch} \, \Hat{Q}_{B_i}^{S_i}(t^{-1})
=(-1)^is_{(1^i)}(\mathrm{\bold x}).
\end{align}
Substituting \eqref{eq-hus} and \eqref{eq-qbs} into \eqref{eq-qus}, we obtain
\begin{align*}
\mathrm{ch} \, \Hat{Q}^{S_{m+d}}_{U_{m,d}}(t)&- t^d \mathrm{ch} \, \Hat{Q}^{S_{m+d}}_{U_{m,d}}(t^{-1})=\sum_{i=0}^{d-1} (-1)^{d}t^i s_{(1^i)}(\mathrm{\bold x})s_{(m+1,1^{d-i-1})}(\mathrm{\bold x})\\
&+\sum_{i=0}^{d-1} \sum_{j=0}^{d-i-1}(-1)^{i+j} t^{d-j} s_{(1^i)}(\mathrm{\bold x}) \Big(s_{(m+d-i-j,1^j)}(\mathrm{\bold x})+s_{(m+d-i-j+1,1^{j-1})}(\mathrm{\bold x})\Big).
\end{align*}
Substituting  $d-j$  for $j$  in the second part of the above equation and then interchanging the order of the summation, we obtain
\begin{align*}
\mathrm{ch} \, \Hat{Q}^{S_{m+d}}_{U_{m,d}}(t)&- t^d \mathrm{ch} \, \Hat{Q}^{S_{m+d}}_{U_{m,d}}(t^{-1})=\sum_{i=0}^{d-1} (-1)^{d}t^i s_{(1^i)}(\mathrm{\bold x})s_{(m+1,1^{d-i-1})}(\mathrm{\bold x})\\
&+\sum_{i=0}^{d-1} \sum_{j=i+1}^{d}(-1)^{i+d-j}t^{j}s_{(1^i)}(\mathrm{\bold x})\Big(s_{(m-i+j,1^{d-j})}(\mathrm{\bold x})+s_{(m-i+j+1,1^{d-j-1})}(\mathrm{\bold x})\Big)\\
&=\sum_{j=0}^{d-1} (-1)^{d}t^j s_{(1^j)}(\mathrm{\bold x})s_{(m+1,1^{d-j-1})}(\mathrm{\bold x})\\
&+\sum_{j=1}^{d} \sum_{i=0}^{j-1}(-1)^{i+d-j}t^{j}s_{(1^i)}(\mathrm{\bold x})\Big(s_{(m-i+j,1^{d-j})}(\mathrm{\bold x})+s_{(m-i+j+1,1^{d-j-1})}(\mathrm{\bold x})\Big).
	\end{align*}
Note that the degree of $\mathrm{ch} \, \Hat{Q}^{S_{m+d}}_{U_{m,d}}(t)$ is strictly less than $\frac{d}{2}$ and hence the degree of lowest term
in $t^d \mathrm{ch} \, \Hat{Q}^{S_{m+d}}_{U_{m,d}}(t^{-1})$  is strict greater than $\frac{d}{2}$.
Thus,
\begin{align}\label{eq-deg0}
[t^0]\mathrm{ch} \, \Hat{Q}^{S_{m+d}}_{U_{m,d}}(t)=(-1)^ds_{(m+1,1^{d-1})}(\mathrm{\bold x}).
\end{align}
Also, for
$1\leq j< d/2$, the coefficient of $t^j$ in $\mathrm{ch} \, \Hat{Q}^{S_{m+d}}_{U_{m,d}}(t)$ is
\begin{align}\label{eq-uni-inv-13}
	(-1)^{d}s_{(1^j)}(\mathrm{\bold x})s_{(m+1,1^{d-j-1})}(\mathrm{\bold x})+\sum_{i=0}^{j-1}(-1)^{i+d-j}s_{(1^i)}(\mathrm{\bold x})\Big(s_{(m-i+j,1^{d-j})}(\mathrm{\bold x})+s_{(m-i+j+1,1^{d-j-1})}(\mathrm{\bold x})\Big).
\end{align}
For any $1\leq j< \frac{d}{2}$, let
\begin{align*}
A(j)=&\sum_{i=0}^{j-1}(-1)^{i+d-j}s_{(1^i)}(\mathrm{\bold x})\Big(s_{(m-i+j,1^{d-j})}(\mathrm{\bold x})+s_{(m-i+j+1,1^{d-j-1})}(\mathrm{\bold x})\Big)\\
	=&\sum_{i=0}^{j-1}(-1)^{i+d-j}s_{(1^i)}(\mathrm{\bold x})s_{(m-i+j,1^{d-j})}(\mathrm{\bold x})+\sum_{i=0}^{j-1}(-1)^{i+d-j}s_{(1^i)}(\mathrm{\bold x})s_{(m-i+j+1,1^{d-j-1})}(\mathrm{\bold x})\\
	=&\sum_{i=0}^{j-1}(-1)^{i+d-j}s_{(1^i)}(\mathrm{\bold x})s_{(m-i+j,1^{d-j})}(\mathrm{\bold x})-\sum_{i=-1}^{j-2}(-1)^{i+d-j}s_{(1^{i+1})}(\mathrm{\bold x})s_{(m-i+j,1^{d-j-1})}(\mathrm{\bold x})\\
=&\sum_{i=0}^{j-2}(-1)^{i+d-j}\Big(s_{(1^i)}(\mathrm{\bold x})s_{(m-i+j,1^{d-j})}(\mathrm{\bold x})-s_{(1^{i+1})}(\mathrm{\bold x})s_{(m-i+j,1^{d-j-1})}(\mathrm{\bold x})\Big)\\
	&~~~~~~~~~~+(-1)^{d-1}s_{(1^{j-1})}(\mathrm{\bold x})s_{(m+1,1^{d-j})}(\mathrm{\bold x})+(-1)^{d-j}s_{(m+1+j,1^{d-j-1})}(\mathrm{\bold x}).
\end{align*}
For $m\geq 1$, $1\leq j< \frac{d}{2}$ and $0 \leq i \leq j-2$, it is easy to verify that $m-i+j\geq 2$  and $(d-j)-i>2$.
By applying Lemma \ref{schur-piere-sub-2}, we obtain
\begin{align*}
s_{(1^i)}(\mathrm{\bold x})s_{(m-i+j,1^{d-j})}(\mathrm{\bold x})-s_{(1^{i+1})}(\mathrm{\bold x})&s_{(m-i+j,1^{d-j-1})}(\mathrm{\bold x})\\
&=-s_{(m-i+j+1,2^i,1^{d-j-i-1})}(\mathrm{\bold x})-s_{(m-i+j,2^{i+1},1^{d-j-i-2})}(\mathrm{\bold x}),
\end{align*}
from which it follows that
\begin{align*}
A(j)=&\sum_{i=0}^{j-2}(-1)^{i+d-j}\Big(-s_{(m-i+j+1,2^i,1^{d-j-i-1})}(\mathrm{\bold x})-s_{(m-i+j,2^{i+1},1^{d-j-i-2})}(\mathrm{\bold x})\Big)\\
	&~~~~~~~~~~+(-1)^{d-1}s_{(1^{j-1})}(\mathrm{\bold x})s_{(m+1,1^{d-j})}(\mathrm{\bold x})+(-1)^{d-j}s_{(m+1+j,1^{d-j-1})}(\mathrm{\bold x})\\
	=&-\sum_{i=0}^{j-2}(-1)^{i+d-j}s_{(m-i+j+1,2^i,1^{d-j-i-1})}(\mathrm{\bold x})
	+\sum_{i=1}^{j-1}(-1)^{i+d-j}s_{(m-i+j+1,2^{i},1^{d-j-i-1})}(\mathrm{\bold x})\\
	&~~~~~~~~~~+(-1)^{d-1}s_{(1^{j-1})}s_{(m+1,1^{d-j})}(\mathrm{\bold x})+(-1)^{d-j}s_{(m+1+j,1^{d-j-1})}(\mathrm{\bold x})\\
	&=-(-1)^{d-j}s_{(m+j+1,1^{d-j-1})}(\mathrm{\bold x})+(-1)^{d-1}s_{(m+2,2^{j-1},1^{d-2j})}(\mathrm{\bold x})\\
	&~~~~~~~~~~+(-1)^{d-1}s_{(1^{j-1})}(\mathrm{\bold x})s_{(m+1,1^{d-j})}(\mathrm{\bold x})+(-1)^{d-j}s_{(m+1+j,1^{d-j-1})}(\mathrm{\bold x})\\
	&=(-1)^{d-1}s_{(m+2,2^{j-1},1^{d-2j})}(\mathrm{\bold x})+(-1)^{d-1}s_{(1^{j-1})}(\mathrm{\bold x})s_{(m+1,1^{d-j})}(\mathrm{\bold x}).
\end{align*}
Therefore, in view of \eqref{eq-uni-inv-13}, for $1\leq j< d/2$ the coefficient of $t^j$ in $\mathrm{ch} \, \Hat{Q}^{S_{m+d}}_{U_{m,d}}(t)$ is
\begin{align}
&(-1)^{d-1}s_{(m+2,2^{j-1},1^{d-2j})}(\mathrm{\bold x})+(-1)^{d-1}s_{(1^{j-1})}(\mathrm{\bold x})s_{(m+1,1^{d-j})}(\mathrm{\bold x})  
+(-1)^{d}s_{(1^j)}(\mathrm{\bold x})s_{(m+1,1^{d-j-1})}(\mathrm{\bold x})\nonumber\\
&=(-1)^{d-1}s_{(m+2,2^{j-1},1^{d-2j})}(\mathrm{\bold x})+(-1)^{d-1} \Big( s_{(1^{j-1})}(\mathrm{\bold x})s_{(m+1,1^{d-j})}(\mathrm{\bold x})-s_{(1^j)}(\mathrm{\bold x})s_{(m+1,1^{d-j-1})}(\mathrm{\bold x}) \Big)\nonumber\\
&=(-1)^{d-1}s_{(m+2,2^{j-1},1^{d-2j})}(\mathrm{\bold x})+(-1)^{d-1} \Big(
-s_{(m+2,2^{j-1},1^{d-2j})}(\mathrm{\bold x})-s_{(m+1,2^{j},1^{d-2j-1})}(\mathrm{\bold x})\Big)\nonumber\\
&=(-1)^ds_{(m+1,2^{j},1^{d-2j-1})}(\mathrm{\bold x}),\label{eq-degj}
\end{align}
where the second last equality is obtained by applying Lemma \ref{schur-piere-sub-2} with the fact that
$m+1\geq 2$ and $(d-j)-(j-1) \geq 2$.
Combining \eqref{eq-deg0} and \eqref{eq-degj}, we obtain
$$[t^j] \mathrm{ch}~\Hat{Q}^{S_{m+d}}_{U_{m,d}}(t)=(-1)^ds_{(m+1,2^{j},1^{d-2j-1})}(\mathrm{\bold x}), \mbox{ for } 0\leq j <d/2,$$
as desired.
This completes the proof. \qed


\subsection{The equivariant Kazhdan-Lusztig polynomials}\label{newproof}

Note that the proof of Theorem \ref{main-thm-inver-uniform} only relies on the evaluation of the equivariant characteristic polynomials for uniform matroids and the inverse Kazhdan-Lusztig polynomials for Boolean matroids. In fact, Theorem \ref{main-eq-o} can be proved in the same manner, assuming that the following result has been proved.

\begin{lem}\label{thm-Boolean-KL}
For any equivariant Boolean matroid $S_n\curvearrowright B_n$  with $n\geq 0$, we have
	$$P_{B_n}^{S_n}(t)=V_{(n)}.$$
\end{lem}
We would like to point out that Lemma \ref{thm-Boolean-KL} is a special case of Theorem \ref{main-eq-o}. Since this lemma can be proved
following the lines of the proof of Theorem \ref{IKL_B_theo}, we omit its proof here. We proceed to prove Theorem \ref{main-eq-o}.

\noindent \textit{Proof of Theorem \ref{main-eq-o}.}
Applying  \eqref{KL-inverse-KL} to the equivariant uniform matroid $S_{m+d} \curvearrowright U_{m,d}$, we obtain
$$\sum_{[F]\in L(U_{m,d})/S_{m+d}} \mathrm{Ind}_{W_F}^{S_{m+d}}
\Big(   P_{(U_{m,d})_{F}}^{(S_{m+d})_F}(t)     \otimes  \Hat{Q}_{(U_{m,d})^F}^{(S_{m+d})_F}(t)   \Big)=0.$$
Recalling the previous arguments on the equivalence of flats in the proof of Theorem \ref{main-thm-inver-uniform} and applying the Frobenius characteristic map lead to
\begin{align}\label{eq-klpuniform}
\sum_{i=0}^{d-1}\mathrm{ch} \, P^{S_{i}}_{B_{i}}(t) \cdot \mathrm{ch} \, \Hat{Q}_{U_{m,d-i}}^{S_{m+d-i}}(t)+ \mathrm{ch} \,P_{U_{m,d}}^{S_{m+d}}(t)=0.
\end{align}
Lemma \ref{thm-Boolean-KL} tells us that, for $0 \leq i \leq d-1$,
	$$\mathrm{ch} \, P_{B_i}^{S_i}(t)=h_i(\mathrm{\bold x})=s_{(i)}(\mathrm{\bold x}).$$
Meanwhile, from Theorem \ref{main-thm-inver-uniform} it follows that
	$$\mathrm{ch} \,  \Hat{Q}_{U_{m,d-i}}^{S_{m+d-i}}(t)=(-1)^{d-i}\sum_{j=0}^{\lfloor (d-i-1)/2 \rfloor}s_{(m+1,2^{j},1^{d-i-2j-1})}(\mathrm{\bold x})t^j.$$
In view of \eqref{eq-klpuniform}, we get
\begin{align*}
 	\mathrm{ch} \,P_{U_{m,d}}^{S_{m+d}}(t)
&=-\sum_{i=0}^{d-1}s_{(i)}(\mathrm{\bold x}) \cdot (-1)^{d-i}\sum_{j=0}^{\lfloor (d-i-1)/2 \rfloor}s_{(m+1,2^{j},1^{d-i-2j-1})}(\mathrm{\bold x})t^j\\
&=-\sum_{i=0}^{d-1}\sum_{j=0}^{\lfloor (d-i-1)/2 \rfloor}(-1)^{d-i}s_{(i)}(\mathrm{\bold x})s_{(m+1,2^{j},1^{d-i-2j-1})}(\mathrm{\bold x})t^j\\
&=-\sum_{j=0}^{\lfloor (d-1)/2 \rfloor}\sum_{i=0}^{d-1-2j}(-1)^{d-i}s_{(i)}(\mathrm{\bold x})s_{(m+1,2^{j},1^{d-i-2j-1})}(\mathrm{\bold x})t^j,
\end{align*}	
where the last equality is obtained by interchanging the order of the summation.

Hence, for $0 \leq j \leq \lfloor (d-1)/2 \rfloor $ the coefficient of $t^j$ in  $ \mathrm{ch} \,P_{U_{m,d}}^{S_{m+d}}(t)$ is
\begin{align}\label{uniform_eq_KL_schur}
[t^j] \mathrm{ch} \,P_{U_{m+d}}^{S_{m,d}}(t)=(-1)^{d+1}\sum_{i=0}^{d-1-2j}(-1)^{i}s_{(i)}(\mathrm{\bold x})s_{(m+1,2^{j},1^{d-i-2j-1})}(\mathrm{\bold x}).
\end{align}
To simplify the summation on the right hand side of \eqref{uniform_eq_KL_schur}, we first consider the constant term of $ \mathrm{ch} \,P_{U_{m,d}}^{S_{m+d}}(t)$, for which we have
\begin{align*}
[t^0] \mathrm{ch} \,P_{U_{m,d}}^{S_{m+d}}(t)&=(-1)^{d+1}\sum_{i=0}^{d-1}(-1)^{i}s_{(i)}(\mathrm{\bold x})s_{(m+1,1^{d-i-1})}(\mathrm{\bold x})\\
&=(-1)^{d+1}  \cdot (-1)^{d-1}s_{(m+d)}(\mathrm{\bold x})\\
&=s_{(m+d)}(\mathrm{\bold x}),
\end{align*}
where the second last equality is obtained by letting $n=d-1$ in Lemma \ref{pieri-4-lem-equi-KL-uniform}.
For $1 \leq j \leq \lfloor (d-1)/2 \rfloor $, letting $n=d-1-2j$ and $m=m+1$ in Lemma \ref{schur-pieri-3-lem}, we find that
\begin{align*}
[t^j] \mathrm{ch} \,P_{U_{m,d}}^{S_{m+d}}(t)
=&(-1)^{d+1} \cdot (-1)^{d-1-2j}\sum_{x=0}^{\mathrm{min}\, \{m-1,d-1-2j \}}s_{(m+1+d-1-2j-x,2+x,2^{j-1})}(\mathrm{\bold x})\\[5pt]
=&\sum_{x=0}^{\mathrm{min}\, \{m-1,d-1-2j \}}s_{(m+d-2j-x,2+x,2^{j-1})}(\mathrm{\bold x})\\[5pt]
=&\sum_{x=1}^{\mathrm{min}\, \{m,d-2j \}}s_{(m+d-2j-x+1,x+1,2^{j-1})}(\mathrm{\bold x}).
\end{align*}
To summarize, we get the desired \eqref{eq-ekl-uniform}. The proof is complete.
\qed

\subsection{A refinement of Lee, Nasr and Radcliffe's formula} \label{newform}

Lee, Nasr and Radcliffe \cite{LNR1-2019, LNR2-2020} considered the combinatorial interpretation for the Kazhdan-Lusztig
polynomials of $\rho$-removed uniform matroids and sparse paving matroids. In particular, they gave the following combinatorial interpretation for the Kazhdan-Lusztig polynomials $P_{U_{m,d}}(t)$ for uniform matroids $U_{m,d}$.

\begin{cor}[{\cite[Theorem 2]{LNR1-2019}}]
For any $m,d \geq 1$, suppose
\begin{align*}
P_{U_{m,d}}(t)=
\sum_{j=0}^{\lfloor (d-1)/2 \rfloor}c_j t^j.
\end{align*}
Then for each $0\leq j\leq \lfloor (d-1)/2 \rfloor$ the coefficient $c_j$ equals the number of standard Young tableaux of shape
$\big(m+d-2j,(d-2j+1)^j\big) \big/ \big( (d-2j-1)^j \big)$.
\end{cor}

We would like to point out that Lee, Nasr and Radcliffe actually used the conjugate partition of $\big(m+d-2j,(d-2j+1)^j\big) \big/ \big( (d-2j-1)^j \big)$ to express $c_j$. Next we will show that the equivariant Kazhdan-Lusztig polynomial $P_{U_{m,d}}^{S_{m+d}}(t)$ admits a more compact form, which implies that the above presentation is more natural in some sense.

\begin{thm}\label{main-equi-schur-skew}
For any $m,d \geq 1$, we have
	\begin{align*}
P_{U_{m,d}}^{S_{m+d}}(t)=
\sum_{j=0}^{\lfloor (d-1)/2 \rfloor}V_{\big(m+d-2j,(d-2j+1)^j\big) \big/ \big( (d-2j-1)^j\big)} t^j.
\end{align*}
\end{thm}
\proof
By Theorem \ref{main-eq-o}, it suffices to show that, for any $m,d,j \geq 1$, we have
\begin{align*}
s_{\big(m+d-2j,(d-2j+1)^j\big) \big/ \big( (d-2j-1)^j \big)}(\mathrm{\bold x})
=\sum_{x=1}^{\mathrm{min}\, \{m,d-2j \}}s_{(m+d-2j-x+1,x+1,2^{j-1} )}(\mathrm{\bold x}).
\end{align*}

Note that, for $\lambda=\big(m+d-2j,(d-2j+1)^j\big)$ and $\mu=\big((d-2j-1)^j \big)$, since the number of cells in the $i$-th row of $\lambda/\mu$ is
$$ \lambda_i-\mu_i=\left\{
\begin{array}{rcl}
m+1,       &      & \mbox{ for }{i=1};\\
2,    &      &\mbox{ for } 2  \leq i \leq j;\\
d-2j+1,       &      & \mbox{ for } {i=j+1},
\end{array} \right.$$
there is a subdiagram of a straight shape $\nu=(m+1,2^j)$ such that all other cells of $\lambda/\mu$ are to the left of $\nu$ and form a single row partition $\rho$.
For example, for $d=10$ and $m=j=3$, the partition $\nu$ is composed of the cells occupied by the bold numbers in Figure \ref{aaaa}, and the partition $\rho$ is composed of other left cells.

\begin{figure}[htb]
\begin{center}
\begin{ytableau}
  \none& \none&  \none&\bold{1}& \bold{1}& \bold{1}& \bold{1} \\
    \none&   \none& \none& \bold{2}& \bold{2} \\
       \none&   \none& \none& \bold{3}& \bold{3} \\
  1&2&2&\bold{4}&\bold{4}
\end{ytableau}
\end{center}
\caption{The Young diagrams of $\lambda/\mu$, $\nu$ and $\rho$}
\label{aaaa}
\end{figure}
In order to get a Littlewood-Richardson tableau $T$ of shape $\lambda/\mu$, there is only one way to fill the cells of
$\nu$ with positive integers, namely, the $i$-th row of $\nu$ is filled with $i$'s.
Moreover, the cells of $\rho$ can only be filled with $1$'s and $2$'s.
Suppose that the number of $2$'s filled in $\rho$ is $x-1$, and hence the number of $1$'s is $d-2j-x$.
According to the lattice permutation condition of $T$, we have $1\leq x\leq d-2j$ and $m+1\geq x-1+2=x+1$,
from which it follows that $1\leq x\leq \mathrm{min}\, \{m,d-2j \}$.
This means that the type $\tau_x$ of $T$ can only be $(m+d-2j-x+1,x+1,2^{j-1})$ for some $1\leq x\leq \mathrm{min}\, \{m,d-2j \}$, and for such a fixed $x$ there is only one Littlewood-Richardson tableau of shape $\lambda/\mu$ and type $\tau_x$. This completes the proof. \qed




\noindent{\bf Acknowledgements.} 
Alice L.L. Gao is supported by  the National Science Foundation of China (No.11801447)  and the Natural Science Foundation of Shaanxi Province (No.2020JQ-104).
Matthew H.Y. Xie is supported by  the National Science Foundation of China (No.11901431).
Arthur L.B. Yang is supported in part by the Fundamental Research Funds for the Central Universities and the National Science Foundation of China (Nos.11522110, 11971249).

\end{document}